\documentclass{article}[12pt]
\usepackage[ruled,vlined]{algorithm2e}
\usepackage{graphicx}
\usepackage{amsmath,amsfonts,amssymb,latexsym}
\usepackage{enumerate}
\usepackage{setspace}
\usepackage{color}
\usepackage[cp1250]{inputenc}
\usepackage{tikz}
\newtheorem{thm}{Theorem}
\newtheorem{prop}[thm]{Proposition}

\newtheorem{obs}[thm]{Observation}

\newtheorem{cor}[thm]{Corollary}

\newtheorem{lem}[thm]{Lemma}

\newtheorem{prob}{Problem}
\newtheorem{ques}{Question}
\newcommand{\qed}{\hfill$\Box$}

\newcommand{\proof}{\noindent\textbf{Proof. }}

\newcommand{\cp}{\,\square\,}

\newcommand{\leaves}{{\mathcal L}}
\newcommand{\supports}{{\mathfrak S}}
\newcommand{\dom}[1]{ {\rm dom}(#1) }
\newcommand{\domg}[1]{ {\rm dom}_g(#1) }
\newcommand{\domga}[1]{ {\rm dom}'_g(#1) }

\begin{document}

\title{The domatic number game played on graphs}
\author{
Bert L. Hartnell $^{a,}$\thanks{Email: \texttt{bert.hartnell@smu.ca}}
\and
Douglas F. Rall $^{b,}$\thanks{Email: \texttt{doug.rall@furman.edu}}
}


\maketitle

\begin{center}
	$^a$  Department of Mathematics and Computing Science, Saint Mary's University,
Halifax, Nova Scotia, Canada\\
	$^b$ Emeritus Professor of Mathematics, Furman University, Greenville, SC, USA	 \\
\end{center}
 \medskip
 
\begin{abstract}
The domatic number of a graph is the maximum number of pairwise disjoint dominating sets admitted by the graph.  We 
introduce a game based around this graph invariant.  The domatic number game is played on a graph $G$  by two players,
Alice and Bob, who take turns selecting a vertex and placing it into one of $k$ sets.  Alice is trying to make each of these sets into a dominating set of $G$ while Bob's goal is to prevent this from being accomplished.    The maximum $k$ for 
which Alice can achieve her goal when both players are playing optimal strategies, is called the game domatic number of $G$.  There are two versions of the game and two resulting invariants depending on whether Alice or Bob is the first to play.

We prove several upper bounds on these game domatic numbers of arbitrary graphs and find the exact values for several classes of graphs including trees, complete bipartite graphs, cycles and some narrow grid graphs.   We pose several open problems concerning the effect of standard graph operations on the game domatic number as well as a vexing question related to the monotonicity of 
the number of sets available to Alice.
\end{abstract}

\noindent
{\bf Keywords:}  domatic number, domination.\\

\noindent
{\bf AMS Subj.\ Class.\ (2020)}: 05C57, 05C69

\section{Introduction}

The \emph{coloring game} was first described by Martin Gardner in his column, \emph{Mathematical Games}~\cite{g-1981} and
then independently defined as a game on graphs by Bodlaender~\cite{b-1991}.  The general game is played on a graph $G=(V,E)$ by two players, Alice and Bob, who take turns coloring an uncolored vertex with a color chosen from a 
set $[k]=\{1,2,\ldots,k\}$ of colors.  When a player assigns a color $c$ to a vertex $x$ no vertex in the neighborhood of 
$x$ can already be colored $c$.  Alice's goal is to eventually color all the vertices of $G$, while Bob's goal is 
to prevent this from happening.  That is, Bob strives to create a situation in which some vertex cannot be chosen and colored 
using any of the $k$ colors.  The \emph{game chromatic number} of $G$, denoted $\chi_g(G)$, is the smallest positive 
integer $k$ for which Alice can achieve her goal regardless of the strategy of Bob. The game chromatic number has been 
investigated by a number of researchers.  For just a small sample see~\cite{bgkz-2007,bfs-2008,dz-1999,h-2025} and their
references.  Recently, a number of games played on graphs that involve some aspect of graph domination have been defined and studied. The outcome of these games typically define a graph invariant.  See~\cite{gik-2019,dgpr-2020,bhkr-2021}.

In this paper we define and study a new game, the \emph{domatic number game}, played on a graph $G$ that is related to the ordinary domatic number.  As in the coloring game, two players, Alice and Bob, alternate turns selecting a vertex from $G$ that has not been selected in the prior moves of the game and assigning a number (shortly called a color) to that vertex.  Alice is attempting to 
assign  colors from a given palette, $[k]$, so that when the game ends every color is present in the closed neighborhood
of every vertex of $G$.  That is, Alice's goal is that for each color $c \in [k]$, the vertices colored $c$ form a dominating set of $G$.  Bob is trying to prevent this. The maximum value of $k$ for which Alice can achieve her goal, regardless of Bob's strategy, is
the \emph{game domatic number} of $G$.  It is denoted $\domg{G}$ when Alice has the first move and by $\domga{G}$ when Bob
has the first move. 

The remainder of the paper is organized as follows.  In the next section we list the most important definitions that will be used 
throughout.  The formal definition of the domatic game and results for general graphs are presented in Section~\ref{sec:gamegeneral}.
Section~\ref{sec:bounds} is devoted to establishing upper bounds for the game domatic number in terms of minimum degree. In Section~\ref{sec:trees} we determine the game domatic numbers for trees, paths, cycles, and some grids.  The paper concludes 
with a number of problems and questions.

\section{Definitions and Notation} \label{sec:defs}
In general, we follow that terminology and notation found in the recent book by Haynes, Hedetniemi and Henning~\cite{hhh-2023}.
If $G$ is a finite, simple graph and $x$ is a vertex of $G$, then $N_G(x)$, the \emph{open neighborhood} of $x$ in $G$, is the set of vertices adjacent to $x$ in $G$.  The \emph{closed neighborhood} of $x$ in $G$
is the set $N_G[x]$ defined by $N_G[x]=N_G(x) \cup \{x\}$.  When the graph is clear from the context we will drop the subscript on these neighborhood names.
A subset $D$ of $V(G)$ is a \emph{dominating set} of $G$ if  $N(x) \cap D \neq \emptyset$, for every $x \in V(G) - D$. The
\emph{domination number} of $G$ is the cardinality of the smallest dominating set of $G$ and is denoted by $\gamma(G)$. For
a positive integer $n$, we denote the set $\{1,\ldots,n\}$ by $[n]$.

A collection of pairwise disjoint subsets of a non-empty set $A$ whose union is $A$ is a \emph{weak partition} of $A$.  In order to simplify terminology, in this paper we will simply call this a \emph{partition} of $A$ even if some of the subsets are empty.
A \emph{domatic partition} of $G$ is a partition, $\Pi$, of $V(G)$ into dominating sets.  If there are $k$ subsets in the partition, then
$\Pi$ is called a $k$-domatic partition.  The maximum $k$ for which $G$ admits a $k$-domatic partition is called the \emph{domatic number} of $G$, and is denoted by $\dom{G}$~\cite{ch-1977}.  Although somewhat dated, see the survey by Zelinka~\cite{z-1998}.
Equivalently, $\dom{G}$ is the maximum number of  pairwise disjoint dominating sets of $G$.
For a given vertex $x$ of $G$ and a dominating set $D$ of $G$, we have $D \cap N[x] \neq \emptyset$.  That is, either $x$ or at least one of its neighbors belongs to $D$.  Applying this to a vertex of minimum
degree in $G$, it follows that $\dom{G} \leq \delta(G)+1$.  A graph $G$ is called \emph{domatically full} if $\dom{G}=\delta(G)+1$.
For a graph $G$ we let $\leaves(G)=\{v \in V(G)\,:\, \deg(v)=1\}$ and let $\supports(G)=\{w \in V(G)\,:\, N(w) \cap \leaves(G) \neq \emptyset\}$.
A vertex in $\leaves(G)$ is a \emph{leaf},  and a vertex in $\supports(G)$ is a \emph{support vertex}.  For $w\in \supports(G)$ we let
$\leaves(w)=N(w) \cap \leaves(G)$.  If a support vertex $w$ is adjacent to more than one leaf, we say it is a \emph{strong support}.
Let $\ell(G)=|\leaves(G)|$ and $\sigma(G)=|\supports(G)|$.

The following result of Ore~\cite{o-1962}, shows that every graph with no isolated vertices has domatic number at least $2$.
\begin{thm} {\rm \cite{o-1962}} \label{thm:Ore}
If $G$ is a graph with no vertices of degree $0$ and $D$ is a minimal dominating set of $G$, then $V(G) - D$ is also a dominating set.
\end{thm}

The \emph{Cartesian product} of two graphs $G$ and $H$ is the graph, $G \cp H$, whose vertex set is $V(G) \times V(H)$.  Two vertices of 
$G \cp H$ are adjacent if they are equal in one coordinate and adjacent in the other.
 The \emph{corona} of $G$ is the
graph $G \circ K_1$ constructed from $G$ by adding a (new) vertex of degree $1$ adjacent to each vertex of $G$.
The \emph{subdivision graph} of $G$, denoted by $S(G)$, is constructed from $G$ by deleting each edge $uv\in E(G)$ and replacing it by a new vertex, $[uv]$, which
is adjacent to both $u$ and $v$.  For $S(G)$ we will refer to the set $V(G)$ as the \emph{old vertices} and $\{ [uv]\,:\, uv \in E(G)\}$ as the \emph{new vertices}.

\section{The Domatic Number Game and General Results} \label{sec:gamegeneral}

Consider the following game played on graphs.  Let $G=(V,E)$ be a graph and let $k$ be a positive integer.  Two players, Alice and Bob, take turns selecting a
previously unchosen vertex $v$ and assigning it (shortly, coloring it with) a color  from $[k]=\{1,\ldots,k\}$.  (Here $[k]$ is called the \emph{palette}.)
We call each such action by a player a \emph{move} or a \emph{play}.   As the game progresses,  we let $V_i$ be the set of vertices that have been colored $i$, for each $i \in [k]$.   When the game has ended,
we call $(V_1,\ldots,V_k)$ the \emph{game induced partition} of $V(G)$.  Alice wins if $(V_1,\ldots,V_k)$ is a domatic partition.  Otherwise, Bob wins.
That is, Alice wins with palette $[k]$ if for each $c \in [k]$ and each $x \in V(G)$ she manages during the course of the game to color at least one vertex from $N[x]$ with color $c$.  Bob wins with palette $[k]$ if at the end of the game there is a vertex $w\in V(G)$ such that at most $k-1$ colors have been assigned to the vertices in $N[w]$.

For a given palette $[k]$, both players follow an optimal strategy to try to win the game. There are two distinct games depending on which player is the first to select and color a vertex.  If Alice has the first move, then we call this the $A$-game and the largest integer $k$ for which Alice has a winning strategy is the \emph{game domatic number}, which is denoted by $\domg{G}$. When Bob is the first to color a vertex we call this the $B$-game, and the largest integer $k$ for which Alice has a winning strategy is the \emph{delayed game domatic number}, which is denoted by $\domga{G}$.  Alice can always win both games with  palette $[1]$.  Therefore, $\domg{G}$ and $\domga{G}$ are well-defined graphical invariants.   For convenience we adopt the following notation.  Suppose the $A$-game is being played on a graph $G$.  When Alice, on her $i^{\rm th}$ move assigns color $j$ to a vertex $v$, we denote this by $A_i(v)=j$.  Similarly,
if Bob assigns color $t$ to vertex $w$ on his $r^{\rm th}$ move we denote this by $B_r(w)=t$.  If the $B$-game was being played, we use $A'_i(v)=j$ and $B'_r(w)=t$ to denote these assignments respectively.  For $X \subseteq V(G)$, we say that a color $c$ is \emph{represented} in $X$ if at least one vertex in $X$ has been assigned color $c$
when the game ends.

The following is a consequence of the fact that the game induced partition is a domatic partition whenever Alice wins the $A$-game or the $B$-game.
\begin{obs} \label{ob:clear}
If $G$ is any graph, then $\domg{G} \leq \dom{G}$ and $\domga{G} \leq \dom{G}$.
\end{obs}

The following lemma will be used a number of times in what follows.

\begin{lem} \label{lem:pmatching}
If a graph $G$ has a perfect matching, then $\domga{G} \ge 2$.
\end{lem}
\proof
Let $M=\{u_1v_1,u_2v_2,\ldots,u_nv_n\}$ be a perfect matching in $G$.  When playing the $B$-game on $G$ with palette $[2]$, Alice can
follow each of Bob's moves by coloring the vertex incident with the matching edge that contains the vertex Bob colored.  That is, if Bob assigns
color $c$ to vertex $u_i$ (respectively, $v_i$), then Alice colors  $v_i$ (respectively, $u_i$) with color $3-c$.  By following this strategy
Alice ensures that the game induced partition $(V_1,V_2)$ is a domatic partition.  Therefore, $\domga{G} \ge 2$.  \qed

The existence of a perfect matching will not necessarily be of use to Alice when the $A$-game is played.  This will be illustrated in
Proposition~\ref{prop:cycles}.

Suppose the $A$-game or the $B$-game is being played on $G$ with palette $[k]$ for some $k\ge 2$.  Suppose $G$ contains a path $u_1u_2u_3u_4$ for which $\deg(u_2)=2=\deg(u_3)$,
$u_2$ and $u_3$ have both been colored $m$ for some $m \in [k]$, and neither $u_1$ nor $u_4$ has yet been colored at some point in the game.  We call this a \emph{$\star mm\star$
configuration}.  It is easy to see that Bob can win this game regardless of which player has the next move.  Indeed, when it is Bob's next move
at least one of the vertices in $\{u_1,u_4\}$ will not be colored, and Bob can assign color $m$ to it.  Also, if  at some point
in a game neither $u_2$ nor $u_3$ has been colored but $u_1$ and $u_4$ have been assigned different colors, say $c$ and $c'$, then we call this a \emph{$c \star\star c'$ configuration}.
If the palette is $[2]$ and Alice is the first to color one of $u_2$ or $u_3$, then it is clear that Bob can win this game.

Also, if $\delta(G)=1$ and the palette is $[k]$ for some $k \ge 2$, then Alice will lose if she is the first player to color a vertex in $\leaves(G) \cup \supports(G)$.  For suppose Alice
colors a support vertex $a$ or an adjacent leaf $b$ with some color $j$.  By coloring the other vertex in $\{a,b\}$ with color $j$ on his next move, Bob will ensure that
no member of the game induced partition other than $V_j$ is a dominating set.  We also note that if the graph $G$ has a strong support vertex $a$, then Bob can, in both the $A$-game and the $B$-game, ensure that only one set in the game induced partition dominates some leaf in $\leaves(a)$.  This establishes the following results for
any graph that has a leaf.

\begin{prop} \label{prop:leaf-parity-strong}
Let $G$ be a graph with minimum degree $1$.
\begin{itemize}
   \item [(i)] If $G$ has a strong support, then $\domg{G}=1=\domga{G}$.
   \item [(ii)] If $G$ is a graph of even order, then $\domg{G}=1$.
   \item [(iii)] If $G$ is a graph of odd order, then $\domga{G}=1$
\end{itemize}
\end{prop}

Using Lemma~\ref{lem:pmatching} and Proposition~\ref{prop:leaf-parity-strong}(ii), the next result follows immediately.

\begin{prop} \label{prop:corona}
If $G$ is any graph, then $\domg{G \circ K_1}=1$ and $\domga{G \circ K_1}=2$.
\end{prop}

\section{Upper Bounds} \label{sec:bounds}

As noted earlier, $\dom{G} \leq \delta(G)+1$.  Thus, $\delta(G)+1$ is an upper bound for both $\domg{G}$ and $\domga{G}$.
By providing an appropriate strategy for Bob, this upper bound can be improved as follows.

\begin{lem} \label{lem:improvedupper}
If $G$ is a graph with no isolated vertices, then
\begin{equation}
  \domg{G}\leq
  \begin{cases}
   \frac{\delta(G)+3}{2},  &\text{if $\delta(G)$ is odd;}\\
   \\
   \frac{\delta(G)+2}{2},  &\text{if $\delta(G)$ is even.}
  \end{cases}
  \text{ and }
  \domga{G}\leq
  \begin{cases}
   \frac{\delta(G)+3}{2},  &\text{if $\delta(G)$ is odd;}\\
   \\
   \frac{\delta(G)+2}{2},  &\text{if $\delta(G)$ is even.}
  \end{cases}
\end{equation}
\end{lem}
\proof
Let the $A$-game be played on $G$ where $m=\delta(G)$.  We establish the upper bounds, based on the parity of $m$, for $\domg{G}$ by providing a strategy for Bob.
  First we suppose that $m$ is odd, say $m=2k+1$, and the palette is $[k+3]$.
We note that $k+3=\frac{\delta(G)+3}{2} +1$.  Bob's strategy is to choose a color $c$ and
assign $c$ to as many vertices in the closed neighborhood of some vertex of degree $m$ as he can.  If Alice's first move is to color a
vertex in $N[x]$, where $\deg(x)=m$, with some color, say $j$, then Bob
picks $c=j$.  Then, following his strategy, Bob can ensure that at least $\frac{1}{2}|N[x]|+1=k+2$ of the vertices in $N[x]$ are assigned color $c$.  The two
players together will have assigned at most $k+1$ distinct colors to the vertices in $N[x]$.  On the other hand, if
on her first move Alice does not color such a vertex, then Bob picks $c=1$ and a vertex $x$ of degree $m$.  Again, following his given strategy, Bob can ensure that at
least $\frac{1}{2}|N[x]|=k+1$ of the vertices in $N[x]$ are colored $c$.  That is, together Alice and Bob will have used at most $k+2$ colors on
the vertices of $N[x]$.
  Consequently, in both cases when the game has ended, $x$ will
not be dominated by at least one of the sets in the game induced partition $(V_1,V_2,\ldots, V_{k+3})$.  Therefore, $\domg{G}\leq k+2=\frac{\delta(G)+3}{2}$.
When the minimum degree of $G$ is even, say $m=\delta(G)=2k$, and the palette is $[k+2]$,  Bob's strategy as above  shows
that  $\domg{G}\leq k+1=\frac{\delta(G)+2}{2}$.

The proof to verify the stated upper bounds in the $B$-game is similar to the above and is omitted.     \qed

It should be noted that no lower bound for $\domg{G}$ and $\domga{G}$ can be given just in terms of 
the minimum degree of a graph.  This is a consequence of the following result of Zelinka.  See Theorem 12.10
in~\cite{hhh-2023}.
\begin{thm}
No minimum degree is sufficient to guarantee the existence of a partition of the vertex set of a graph into three 
dominating sets.
\end{thm}

Note that if the graph is regular of odd degree, then the upper bound for the game domatic number from Lemma~\ref{lem:improvedupper} can be improved since
Alice's first move will be to color a vertex in the closed neighborhood of a vertex of minimum degree.
\begin{cor} \label{cor:oddregular}
If $G$ is regular of odd degree, then $\domg{G} \leq \frac{\delta(G)+1}{2}$.
\end{cor}

If a graph has more than one vertex that dominates the entire graph, then Alice can, by playing as many of these as possible, ensure that the game
domatic numbers have the lower bound given in the next result.   
\begin{obs}
If a graph $G$ has $k\geq 2$ universal vertices, then $\domg{G} \geq \lceil{k/2}\rceil$ and $\domga{G} \geq \lceil{k/2}\rceil$.
\end{obs}

The following is straightforward to verify and also shows that the bounds in Lemma~\ref{lem:improvedupper} are sharp.  In addition,
the class of complete graphs shows that $\dom{G}-\domg{G}$ can be arbitrarily large.

\begin{prop} \label{prop:cliques}
If $n$ is a positive integer, then
\begin{equation}
  \domg{K_n}= \Bigl \lceil \frac{n}{2}\Bigr \rceil \text{ and } \domga{K_n}=
  \begin{cases}
   \frac{n+1}{2},  &\text{if $n$ is odd;}\\
   \frac{n+2}{2},  &\text{if $n$ is even.}
  \end{cases}
\end{equation}
\end{prop}

\section{Trees and Other Graph Classes} \label{sec:trees}

\begin{prop} \label{prop:grids}
Let $2 \leq m \leq n$.  If at least one of $m$ and $n$ is even, then $\domga{P_n \cp P_m}=2$.  Also, $\domg{P_n \cp P_2}=1$.
\end{prop}
\proof
Since $\delta(P_n \cp P_m)=2$, it follows from Lemma~\ref{lem:improvedupper} that $\domga{P_n \cp P_m} \leq 2$.
Suppose first that at least one of $m$ and $n$ is even.  In this case the grid $P_n \cp P_m$ admits a perfect matching, and hence
by Lemma~\ref{lem:pmatching}, $\domga{P_n \cp P_m}=2$.

Let the $A$-game be played on the grid $P_n \cp P_2$.
For ease of reference in the Cartesian product $P_n \cp P_2$ we let $x_k=(k,1)$ and $y_k=(k,2)$ for $k \in [n]=V(P_n)$.
By Lemma~\ref{lem:improvedupper}, we have $\domg{P_n \cp P_2}\leq 2$.  We provide a strategy for Bob that shows this upper bound is not
attained when the palette is $[2]$.  Suppose on her first move Alice colors one of the
vertices of degree $2$.  Without loss of generality we may assume that $A_1(y_1)=1$.  Bob responds with $B_1(x_1)=1$, which creates a
$\star 11\star$ configuration and thus guarantees that Bob can win the game.
Therefore, if Alice's first move is to color a ``corner'' vertex, then Bob wins the game.  Suppose then that $n \ge 3$ and
Alice's first move is to assign a color, say $1$, to $y_k$ for some $k \in [n]$ with $1<k<n$. Bob's first move then is to let $B_1(x_k)=2$.
Alice avoids letting Bob create a $\star cc \star$ configuration, and thus she will color a vertex of degree $3$ until no such vertex is uncolored.  For
all such moves (that is, for $2 \leq i \leq n-2$), if $A_i(y_j)=c$ (respectively, $A_i(x_j)=c$), then Bob responds with $B_i(x_j)=3-c$ (respectively, $B_i(y_j)=3-c$).  By following this strategy
Bob will, after $2n-4$ moves are made in the game, force Alice to be the first player to color a vertex in a $c\star\star\,c'$ configuration.
By following this strategy Bob can win the game when the palette is $[2]$.
Therefore, $\domg{P_n \cp P_2}=1$.  \qed

Every nontrivial tree is domatically full. That is, if $T$ is a tree of order at least $2$, then $\dom{T}=2$. As the following two results show, there is only a small class of trees whose game domatic number or delayed game domatic number achieves this upper bound.  By Proposition~\ref{prop:leaf-parity-strong},
if a tree $T$ has a strong support vertex, then $\domg{T}=1=\domga{T}$.  Since every nontrivial tree has domatic number $2$,
in every proof involving a tree we assume we are playing the game with $[2]$ as the palette.

\begin{thm} \label{thm:AB-game-trees}
If $T$ be a tree of order at least $2$.   Then
\begin{itemize}
\item [(i)] $\domga{T}=2$ if and only if $T$ has a perfect matching, and
\item [(ii)] $\domg{T}=1$.
\end{itemize}
\end{thm}
\proof
Let the $B$-game be played on $T$.  We may assume that $T$ has no strong support vertices.
If $T$ admits a perfect matching, then by using Lemma~\ref{lem:pmatching} we see that $\domga{T}=2$.  For the converse we assume
that $T$ has no perfect matching.
If $T$ has odd order, then Bob can force Alice to be the first player to color a vertex in $\leaves(T) \cup \supports(T)$.
It follows that Bob wins the game and $\domga{T}=1$.  Thus we may assume that $T$ has even order.
We define a  sequence of forests as follows.
Let $F_0=T$, and for  $i \geq 1$, let $F_i=F_{i-1}-(\leaves(F_{i-1}) \cup \supports(F_{i-1}))$. Repeat this process until some forest $F_q$ has a component
that is either an isolated vertex or a star $K_{1,r}$ for some $r \ge 2$. (We will call such a star ``large''.)   There must exist a smallest such $q$ since $T$ does not have a perfect matching.  Bob plays the following strategy.  He orders the support vertices of $F_0=T$, and in consecutive moves assigns
color $1$ to each.  After each such move of Bob, Alice, not wanting to lose, will be forced to color the adjacent leaf with color $2$.  If $F_1$ has no
component that is an isolated vertex or a large star, then Bob orders the support vertices of $F_1$ and colors them as he did in $F_0$.  Again Alice will be
forced to color the corresponding leaves with color $2$.  This is true since if $x \in \leaves(F_1)$, then the neighbors of $x$ in $T$ belong to $\supports(F_1)\cup \supports(F_0)$.
That is, after Bob colors the support vertex from $F_1$ that is the unique neighbor of $x$ in $F_1$, every vertex in $N_T(x)$ has been colored $1$.
Bob continues this approach until the vertices in $(\supports(T)\cup \leaves(T)) \cup \cdots \cup (\supports(F_{q-1})\cup \leaves(F_{q-1}))$ have been colored, and it is Bob's turn.  If $F_q$ has an isolated vertex $z$,
then Bob assigns it color $1$.  All of the vertices in $N_T[z]$ are now colored $1$.  Otherwise, $F_q$ has a component that is a large star with center $w$.  Bob assigns color $1$ to $w$.  After Alice's next
move there will be at least one uncolored leaf, say $y$, of this large star.  Bob assigns color $1$ to $y$.  When the game ends all the vertices
in $N_T[y]$ are colored $1$.  That is, $V_2$ is not a dominating set of $T$.  Therefore, $\domga{T}=1$.  This finishes the proof of $(i)$.

Now let the $A$-game be played on $T$.   If $T$ has even order, then it follows by Proposition~\ref{prop:leaf-parity-strong} that $\domg{T}=1$.
Hence, we may assume that $T$ has odd order.  On her first move Alice will not color a vertex in $\supports(T)\cup \leaves(T)$ as was observed in the
paragraph preceding Proposition~\ref{prop:leaf-parity-strong}.  Thus, without loss of generality we may
assume that Alice assigns color $1$ to a vertex $x \in V(T)-(\leaves(T) \cup \supports(T))$.  Bob now follows a strategy similar to that outlined in the proof of the
last case in $(i)$ above.  That is, Bob determines the sequence of forests $F_0, F_1, \ldots$ and the smallest index $q$ such that $F_q$ has a component
of order $1$ or a large star.  

Suppose first that $x \notin \cup_{i<q}(\supports(F_i) \cup \leaves(F_i))$. 
For each $i \leq q-1$, Bob assigns color $1$ to the support vertices of $F_i$.  In those forests Alice will assign color $2$
to the corresponding leaves so that she does not at that point in the game create a closed neighborhood contained in $V_1$.  Thus we arrive
at the forest $F_q$.  If there exists an isolated vertex, then it is either $x$, whose closed neighborhood has now been colored $1$ or is an uncolored vertex, say
$z$.  In the latter case, Bob assigns color $1$ to $z$.  In both cases $V_2$ will not be a dominating set when the game ends.  Otherwise, $F_q$ has a large
star with center $w$.  If $w=x$, then Bob assigns color $1$ to any leaf in this star.  If $w \neq x$, then Bob colors $w$ with color $1$.  Since the star has
more than one leaf in $F_q$, it follows that Bob can ensure that the closed neighborhood of one of its leaves will be contained in $V_1$ when the game ends.

Finally assume that the vertex $x$ appears as a leaf or a support vertex in $F_r$ for some $r<q$. As above, Bob assigns color $1$ to each support vertex in $F_i$ for
$0 \le i \le r-1$.  If $x \in \leaves(F_r)$, then Bob assigns color $1$ to its adjacent support vertex in $F_r$. This implies that every vertex in $N_T[x]$ will be colored $1$ at the end of the game. On the other hand, if $x \in \supports(F_r)$ and the unique neighbor of $x$ that is a leaf in $F_r$ is $y$, then Bob assigns color $1$ to $y$.  In this case all vertices in $N_T[y]$ will be colored $1$ at the end of the game. 
In both cases, $V_2$ will not be a dominating set of $T$ when the game has ended.  Therefore, $\domg{T}=1$.   \qed

The next result follows directly from Theorem~\ref{thm:AB-game-trees}.
\begin{cor} \label{cor:paths}
Let $n$ be a positive integer. Then
\begin{equation}
  \domg{P_n}=1 \text{  and  }
  \domga{P_n}=
  \begin{cases}
   1,  &\text{if $n$ is odd;}\\
   2,  &\text{if $n$ is even.}
  \end{cases}
\end{equation}
\end{cor}

\begin{prop} \label{prop:cycles}
Let $n$ be a positive integer with $n \geq 3$.  Then
\begin{equation}
  \domg{C_n}=
  \begin{cases}
   2, &\text{ if $n=3$;} \\
   1, &\text{ if $n \geq 4$.}
  \end{cases}
   \text{  and  }
  \domga{C_n}=
  \begin{cases}
   1,  &\text{if $n\geq 5$ is odd;}\\
   2,  &\text{if $n=3$ or $n$ is even.}
  \end{cases}
\end{equation}
\end{prop}
\proof
Let $V(C_n)=\{v_i \,:\, i \in [n]\}$ and $E(C_n)=\{v_iv_{i+1} \,:\, i \in [n] \text{ subscripts modulo } $n$\}$.  In either game played on a
cycle, if Bob can employ a strategy that ensures three consecutive vertices on the cycle are assigned the same color, then he wins that game.  Clearly, $\domg{C_3}=2$.
Let $n \geq 4$ and let the $A$-game be played on $C_n$ with palette $[2]$.   Without loss of generality we may assume that on her first move Alice assigns color $1$ to $v_1$.
Bob then assigns color $1$ to $v_2$, which creates a $\star11\star$  configuration and Bob wins the game.  Therefore, $\domg{C_n}=1$.

Clearly, $\domga{C_3}=2$.
Now let the $B$-game be played on $C_5$.  We will show that Bob can ensure that three consecutive vertices on $C_5$ are assigned the same color.  We may assume Bob assigns color $1$ to $v_1$ on his first move.  Clearly, Alice will not assign color $1$ to $v_2$ or to $v_5$, either move of which would create a $\star11\star$ configuration.
If, on her first move, Alice colors $v_2$ with color $2$, then Bob can color $v_5$ with color $1$,  which will force Alice to assign color $2$ to $v_4$ to prevent $v_4,v_5$ and $v_1$ from all being assigned color $1$. However, then Bob assigns color $2$
to $v_3$, and Bob wins since $v_3$ is dominated by color $2$ only.    By symmetry
Alice will not color $v_5$ on her first move.  Therefore, we may assume by symmetry that Alice's first move will be to color $v_3$.  She will not color $v_3$
with color $1$ since that would allow Bob to assign color $1$ to $v_2$.  Therefore, she will assign color $2$ to $v_3$.  Bob can then assign color $1$ to $v_5$,
thereby creating a $\star11\star$ configuration, and Bob can win the game.  Therefore, $\domga{C_5}=1$.

Now let $n=2p+1$ with $p \geq 3$ and let the $B$-game with palette $[2]$ be played on $C_n$.  We exhibit
a strategy for Bob to ensure that he wins this game.  Bob's first move is to assign color $1$ to $v_1$.  In order to avoid Bob creating a $\star11\star$ configuration, Alice must assign color $2$ to either $v_2$ or $v_n$.  By symmetry we may assume $A'_1(v_2)=2$.  On his second move Bob plays $B'_2(v_5)=1$,  which creates a $2\star\star1$ configuration on $v_2v_3v_4v_5$.  Since an even number of vertices remain uncolored in the set $\{v_j\,:\, 6 \leq j \leq n\}$ and it is Alice's turn, Bob can force Alice to be the first to color $v_3$ or $v_4$.  Therefore, Bob wins, and $\domga{C_n}=1$.

Now let the $B$-game be played on $C_n$, where $n$ is even.  By Lemma~\ref{lem:improvedupper}, $\domga{C_n}\leq 2$.  By applying Lemma~\ref{lem:pmatching},
we conclude that $\domga{C_n} = 2$ since $C_n$ admits a perfect matching.  \qed

\begin{thm} \label{thm:completebipartite}
If $2 \le m \le n$, then
\begin{equation}
  \domg{K_{m,n}}=
  \begin{cases}
   \frac{m}{2}, &\text{ if $m$ and $n$ are both even;} \\
   \lceil\frac{m+1}{2}\rceil, &\text{ otherwise.}
  \end{cases}\\
\end{equation}
and  
\begin{equation}
  \domga{K_{m,n}}=
  \begin{cases}
   \frac{m}{2},  &\text{if $m$ is even and $n$ is odd;}\\
   \lceil\frac{m+1}{2}\rceil,  &\text{otherwise.}
  \end{cases}
\end{equation}
\end{thm} 

\proof 
Let $V=\{v_1,v_2, \ldots,v_m\}$ and $W=\{w_1,w_2,\ldots,w_n\}$ be the two maximal independent sets of $K_{m,n}$.
For convenience we will refer to $V$ and $W$ as the ``sides.''
Without loss of generality we may assume that the first player to color a vertex in either game uses color $1$.
In describing a player's strategy we will say that the player ``follows'' the other player to mean he or she
colors a vertex from the same side as the other player just did on the immediately preceding move.
The general strategy for Bob is the following.  If he is going to play the vertex $x$, then he assigns color $1$ to $x$
if no vertices from that side have been colored previously.  Otherwise, he uses the smallest color already represented 
in that side.  Alice's general strategy, on every move except when she begins the $A$-game, is the following.  If she is playing a vertex $v_i$ (respectively, $w_i$) and the largest color 
represented in $V$ (respectively, $W$) up to that point in the game is color $j$, then she assigns color $j+1$ to $v_i$ (respectively, $w_i$).

{\bf Case 1: [$m$ and $n$ both even]} Let the $A$-game be played on $K_{m,n}$ with palette $[k]$ where $k=m/2$.
Alice's first move is to assign color $1$ to $w_1$.  For the remainder of the game,  Alice follows Bob  if this is possible.  If it is not possible, then she colors a vertex in the other side.  By doing this 
and employing her general strategy described above, Alice can win the game with palette $[k]$. Suppose now that the palette has more than $k$ colors. Since $m$ and $n$ are both even, Bob can force Alice to be the first player to play a vertex in $V$.  By following his general strategy outlined above, Bob can ensure that at most $k$ colors are represented in $V$ when the game ends.  
Thus Bob wins the $A$-game when the larger palette is used.  It now follows that $\domg{K_{m,n}}=m/2$ when $m$ and $n$ are both even.  

Now let the $B$-game be played on $K_{m,n}$ with palette $[k]$, where $k=\lceil\frac{m+1}{2}\rceil=\frac{m}{2}+1$.  Since $n$ is even, Alice can force Bob to be the first one to play on $V$.  By following Bob and employing her general strategy given above, she can ensure that all colors from the palette are represented in both $V$ and $W$ at the end of the game.  As above, if the palette has more than $\frac{m}{2}+1$ colors, then Bob
can win the game by following his general strategy.  Therefore, $\domga{K_{m,n}}=\lceil\frac{m+1}{2}\rceil$ in this case.

{\bf Case 2: [$m$ even and $n$ odd]} Let the $A$-game be played on $K_{m,n}$ with palette $[k]$, where $k=\lceil\frac{m+1}{2}\rceil=\frac{m}{2}+1$.  Alice's first move is to assign color $1$ to $w_1$.  Then, by playing her general strategy she can force Bob to be the first player to color a vertex in $V$ since $n$ is odd.  Regardless of how Bob plays, Alice can follow Bob and use her general strategy for the remainder of the game, which means that all $k$ colors will be represented in both sides 
when the game ends. By coloring a vertex in $V$ on his first move and then using his general strategy, Bob can prevent more than $k$ colors from being represented in $V$.  Therefore, $\domg{K_{m,n}}=\lceil\frac{m+1}{2}\rceil$ when $m$ is even and $n$ is odd.

Suppose the $B$-game is played on $K_{m,n}$ with palette $[k]$, where $k=\frac{m}{2}+1$.  Bob's first move is to assign color $1$ to $w_1$.  Since $n$ is odd, he can force Alice to be the first player to color a vertex in $V$.  Hence, by using his general strategy
and always following Alice, at the end of the game at most $k-1$ of the colors will be represented in $V$.  It follows that the game induced partition will not be a domatic partition.  Bob wins this game, and this implies that $\domga{K_{m,n}}\leq k-1=\frac{m}{2}$.  With a palette of $[\frac{m}{2}]$ Alice can win the game using the following approach.  After Bob's first move, she makes her first move by using her general strategy and coloring a vertex in $V$.  After this she follows Bob and uses her general strategy.  This will ensure that all $\frac{m}{2}$ colors are represented in both $V$ and $W$.  Therefore, $\domga{K_{m,n}}= \frac{m}{2}$ when $m$ is
even and $n$ is odd.

{\bf Case 3: [$m$  odd]}  Let the $A$ game be played on $K_{m,n}$ with palette $[k]$, where $k=\lceil\frac{m+1}{2}\rceil$.  Alice's first move is to assign color $1$ to $v_1$.  For the remainder of the game she always follows Bob and uses her general strategy.  It is easy to see
that she can ensure all $k$ colors are represented in both sides.  If the palette contains $\ell$ colors, where $\ell \geq k+1$, then Bob can win the game by continually assigning color $1$ to uncolored vertices in $V$ until being forced to play in $W$.  
He then assigns color $1$ to at least one vertex in $W$.  Regardless of how the game proceeds after that, at most $k$ colors will be represented in $V$ and the game induced partition  
$(V_1,\ldots,V_{\ell})$ is not a domatic partition.  That is, $\domg{K_{m,n}}=\lceil\frac{m+1}{2}\rceil$.

Finally, suppose the $B$-game is played on $K_{m,n}$ with $[k]$ as the palette, where $k=\lceil\frac{m+1}{2}\rceil=\frac{m+1}{2}$.  It is easy to see that Alice can win with this palette by always following Bob and using her general strategy.  With a larger palette, Bob can win the
game by using his strategy given for the $A$-game in this case.  Therefore, $\domga{K_{m,n}}=\lceil\frac{m+1}{2}\rceil$.  This completes the proof. \qed

Another class of graphs for which the game domatic invariants can be determined are certain subdivision graphs.  Since they all have minimum
degree $2$, it follows from Lemma~\ref{lem:improvedupper} that we may use $[2]$ as the palette.
\begin{lem} \label{lem:overlappingcycles}
If a graph $G$ contains two edge-disjoint cycles that have a single vertex in common, then $\domga{S(G)}=1$.
\end{lem}
\proof
We define a strategy for Bob when the $B$-game is played on $S(G)$.  Let the two cycles of $G$ be $P=x,v_1,v_2,\ldots,v_p,x$ and $Q=x,u_1,u_2,\ldots,u_q,x$.
Bob's opening move is $B_1'(x)=1$.  Since $x$ is the only vertex shared by the two subdivided cycles, we may assume without loss of generality that on her first move
Alice does not color a vertex on $P$ or any of the new vertices created when the edges of $P$ are subdivided.  Bob's second move is to assign color $1$
to $v_1$. This will require Alice to respond with $A_2'([xv_1])=2$.  Bob's strategy is to then assign color $1$ to $v_2$, which forces Alice to play
$A_3'([v_1v_2])=2$. Bob continues ``around the subdivided $P$'' by assigning each consecutive old vertex the color $1$.  After Bob's $(p+1)^{\rm st}$ move all three of the vertices in the set $\{v_{p-1},v_p,x\}$ have been assigned color $1$, and neither of $[v_{p-1}v_p]$ or $[v_px]$ have been colored.  It is clear that Bob can assign color
$1$ to one of these on his next move.  Therefore, $V_2$ will not dominate at least one these two new vertices when the game ends.  We conclude that $\domga{S(G)}=1$.  \qed

Any grid graphs with each dimension at least $3$ has two edge-disjoint cycles that share a single vertex.
\begin{cor}
If $3 \leq m \leq n$, then $\domga{S(P_m \cp P_n)}=1$.
\end{cor}

A slight modification of the proof of Lemma~\ref{lem:overlappingcycles} can be used to prove the following.  After Alice's first move, Bob can color one of the two common vertices and then proceed with the strategy described for him in the proof of Lemma~\ref{lem:overlappingcycles}.
\begin{lem}
If a graph $G$ has four pairwise edge-disjoint cycles $F_1,F_2,F_3$ and $F_4$ such that $F_1$ and $F_2$ share a single vertex and $F_3$ and $F_4$ share a
single vertex while $F_1 \cup F_2$ and $F_3 \cup F_4$ are vertex disjoint, then $\domg{S(G)}=1$.
\end{lem}
\begin{cor}
If $m \ge 3$ and $n \ge 6$, then $\domg{S(P_m \cp P_n)}=1$.
\end{cor}

\section{Questions and Further Research}
We close with some questions and problems that we did not address in this initial study of
the domatic number game.

Deleting an edge from a graph $G$ can lower the (ordinary) domatic number.  For example, this is true for any nontrivial
complete graph or any cycle of order a multiple of $3$.  The following question is natural for the game domatic number.     
\begin{ques} How does edge or vertex removal  affect the game domatic number?  
\end{ques}

\begin{prob}
Find lower and upper bounds for $\domg{G \cp H}$ and $\domga{G \cp H}$ in terms of graphical invariants of arbitrary graphs  $G$ and $H$.
\end{prob}

One of the curious outstanding questions in the study of the coloring game is the following; see~\cite{h-2025}.  If Alice 
can win the coloring game on a graph $G$ when there are $k$ colors available, can she also win when $k+1$ colors are 
available?    A similar vexing question arises in the domatic number game.
\begin{ques}
Does there exist a graph $H$ and a positive integer $k$ such that Bob wins the $A$-game (respectively, $B$-game) played on $H$ with palette $[k]$ but Alice wins the $A$-game (respectively, $B$-game) played on $H$ with palette $[k+1]$?
\end{ques}

\begin{prob}
Find an expression for $\domg{G_1 \cup G_2}$ and $\domga{G_1 \cup G_2}$ in terms of $\domg{G_i}$
and $\domga{G_i}$ for $i \in [2]$ when $G_1$ and $G_2$ are arbitrary vertex disjoint graphs.
\end{prob}

\end{document}